\documentclass[journal]{IEEEtran}
\ifCLASSINFOpdf
\else
\fi
\usepackage{graphicx}
\usepackage{amsmath}
\usepackage{amsthm,enumerate}
\newtheorem{definition}{Definition}
\newtheorem{theorem}{Theorem}

\newtheorem{lemma}{Lemma}

\usepackage{amssymb}
\usepackage{enumitem}

\begin{document}
%
\title{An Overview of Energy-Optimal Impedance Control of Cooperative Robot Manipulators}
%
%
%

\author{Amin~Ghorbanpour,
        Hanz~Richter
\thanks{The authors are with the Department of Mechanical Engineering, Cleveland
State University, Cleveland, OH 44122 USA (e-mail: a.ghorbanpour@csuohio.edu;
h.richter@csuohio.edu).}
\thanks{This article is intended as a supplement to article titled ``Energy-Optimal Impedance Control of Cooperative Robot Manipulators" written by the same authors.}
}

\maketitle

\begin{abstract}
An impedance-based control scheme is introduced for cooperative manipulators grasping a rigid load. The position and orientation of the load are to be maintained close to a desired trajectory, trading off tracking accuracy by low energy consumption and maintaining stability. To this end, the augmented dynamics of the robots, their actuators and the load is formed, and an impedance control is adopted. A  virtual control strategy is used to decouple torque control from actuator control. An optimization problem is then formulated using energy balance equations. The optimization finds the damping and stiffness gains of the impedance relation such that the energy consumption is minimized. Furthermore,  $\mathcal{L}_{2}$ stability techniques are used to allow for time-varying damping and stiffness in the desired impedance. A numerical example is provided to demonstrate the results.
\end{abstract}


%
\IEEEpeerreviewmaketitle

\section{Introduction}
%
%
%
%
\IEEEPARstart{S}{aving} energy has become a major driving force in engineering, owing to increasing demand for cost-efficient and long-life systems, and as a result, engineers sought to find solutions for optimal energy consumption. This is particularly important in robotic automation. Robots with on-board, finite energy storage are prevalent in electric vehicles, powered human assistive devices, aerospace vehicles, etc. \cite{khalaf2019trajectory,richter2015framework,ghorbanpour2018control,khalaf2016parametric,shimizu2013super,khalaf2018development,richter2020control,he2019research,boscariol2019trajectory,khalaf2018global,richter2015impedance}.

 

Therefore, in dynamical systems with limited on-board energy source, longer operating times and smaller operating costs are possible by reducing the amount of energy needed to complete the motion tasks. 
Energy regeneration is defined here as the ability to transfer energy back to the system's power source during motion tasks.  
Energy-aware control in robotics, particularly concerning energy regeneration, has gained a broad research interest recently~\cite{khalaf2019trajectory,richter2015framework,he2019research,boscariol2019trajectory,khalaf2018development}.
Ultracapacitors have been used in addition to or instead of batteries as energy storage elements in regenerative motion systems. In addition to being lightweight and durable, ultracapacitors can be charged and discharged at high rates. This feature is essential for energy regeneration because the storage element should be capable of admitting all the energy that can be potentially recovered. 

So far, most studies have considered energy-oriented control in single robots where  significant potential for energy recovery exists, e.g. industrial robots, prosthetic legs, etc. \cite{richter2020control,khalaf2019trajectory,khalaf2018development}. These works sought to minimize energy consumption by maximizing regenerative energy storage relative to trajectories, controller parameters or design parameters. The latter problem was quadratic, and admitted explicit solutions. 

The focus of this study is to introduce a framework for energy-oriented control in cooperative robots. Even though robots offer  superior capabilities than humans for carrying out motion tasks involving large loads, speeds and tight precision requirements, some tasks are difficult or impossible for a single robot. The limits in the structure of a robot prevent a sole robot to operate with large, unbalanced or flexible loads. Cooperative robot manipulators (CRM) offer significant advantages over a single robot, and they show better performance in tasks such as grasping, gripping, lifting, transferring, lowering, and approach an object. Regardless of the benefits obtained by employing CRM, it comes
at the cost of complexity for the mathematical modeling, control, and coordination. Here, the term cooperative means collaboration of multi robots in handling a payload.

This paper considers theoretical developments and extensions of a framework previously established for a single robot in \cite{richter2015framework}. Here, a motion task is defined as carrying a payload near a desired trajectory by grasping it using the CRM, limiting grasping forces and reducing energy consumption. The motion task thus demands a controller meeting three main objectives: 1) maintaining acceptable tracking accuracy of desired trajectories for the position and orientation of the load, 2) being in compliance with the closed kinematic chain caused by grasping the load, and 3) minimizing the energy consumption through energy regeneration. A novel  controller is proposed based on an impedance  relation whose parameters are tuned such that maximum energy regeneration is provided. It is assumed that all robots and their corresponding actuators (DC motors) are equipped with regenerative drives (4-quadrant drives), allowing energy to go back from the robots to the storage element. Regenerative drives provide an  opportunity to harvest the excess mechanical energy by channeling it back to the source instead of being dissipated \cite{khalaf2019trajectory}.

\subsection{Review of related work}
Most of the research literature on CRM does not place  emphasis on energy considerations. Optimization of the energy consumption for single robots and general electromechanical systems is well documented \cite{shimizu2013super,richter2015framework,richter2020control,khalaf2019trajectory,richter2015impedance,khalaf2016parametric,khalaf2018development,ghorbanpour2018control,warner2016design,carabin2017review,richter2020optimal,joy2018regenerative,ivoilov2018power}. In the field of robotics, researches in \cite{richter2015framework,khalaf2019trajectory,khalaf2016parametric,ghorbanpour2018control} minimized energy consumption considering capacitive storage elements. These ideas were applied for powered human assistive devices~\cite{khalaf2018development,richter2015impedance,warner2016design}, and aerospace applications~\cite{richter2020control,richter2020optimal,shimizu2013super}. Energy-oriented control studies have also been conducted with a focus on actuator design and dynamics \cite{ghorbanpour2018control,joy2018regenerative,ivoilov2018power}. 

The framework established in~\cite{richter2015framework} combines the dynamics of a robot, joint mechanisms, and actuators to present an augmented dynamical model of the robotic system. Depending on the type of the actuators, e.g. DC motors, hydraulic elements, etc., the energy transfer between the power source and the actuators can be derived, which leads to find the energy change in the power source. Hence, an energy minimization can be formed. These ideas were applied to an industrial robot (PUMA manipulator), finding point-to-point energy-optimal trajectories~\cite{khalaf2019trajectory}. It is shown that 10-20\% energy reduction is possible by energy regeneration. Furthermore, application of the idea for prosthetic leg showed promising results \cite{khalaf2018development,khademi2017optimal,warner2016design}. 

In the aforementioned works, robots are either isolated without any contact with the environment, or have limited (in the sense of low and intermittent interaction forces) interactions with it. Accordingly, interactions were not modeled and regarded as external disturbances. This work provides a bridge between energy-oriented control works and the subject of cooperative robots, considering the stability of interactions.

Modeling and control of CRM have received wide attention for more than 3 decades\cite{koivo1990modeling,chiacchio1992cooperative,walker1991analysis,jean1993adaptive,erhart2015internal,vergnano2012modeling,kawasaki2006decentralized,bonitz1996internal,caccavale2008six,szewczyk2002planning,shimizu2012nonlinear,he2016dual,heck2013internal,schneider1989object,li2017indirect,ren2014adaptive,caccavale2001impedance,erhart2013impedance,stolfi2017combined}. Nevertheless, the energy-oriented studies of CRM have received little attention. For instance, in \cite{vergnano2012modeling} energy consumption is reduced by finding the optimal operating schedule of a robotic manufacturing system, using the consumption profiles of individual robots and their operations as pre-determined information.

CRM tasks involve interaction between the robots and a payload, and the force-position relationships between them are of fundamental concern. For trajectory tracking, various control approaches have been proposed, with impedance control being prevalent. Impedance control constitutes an effective, easy-to-implement scheme for multiple interacting dynamical systems. Impedance control includes regulation and stabilization of robot motion by establishing a mathematical relationship between the interaction forces/moments and the robot trajectories. The relation usually defined as a second order linear non-homogeneous ordinary differential equation, often modeled after a forced spring-mass-damper system. Impedance control has been used for CRM in many works, for instance ~\cite {caccavale2008six,caccavale2001impedance,bonitz1996internal,szewczyk2002planning,li2017indirect,ren2014adaptive,he2016dual,erhart2013impedance,stolfi2017combined,heck2013internal,schneider1989object,shimizu2012nonlinear}.

One of the main differences among various approaches is the definition of terms and the structure of the impedance relationship. This includes different choices of the forces or moments appearing in the impedance relation~\cite{caccavale2008six,bonitz1996internal} or the use of nonlinear~\cite{shimizu2012nonlinear} and variable~\cite{he2016dual} inertial, stiffness and damping characteristics.

Furthermore, impedance control is useful when the payload handled by the CRM may also be in contact with the environment. The algorithm of~\cite{heck2013internal} controls motion and internal forces of the payload, as well as the contact forces between the payload and environment. In that work, the reference trajectory for each robotic joint is calculated based on an impedance relationship such that the desired internal and contact forces are achieved. In~\cite{schneider1989object} the problem of assembly operations using CRM was considered. Experiments showed the effectiveness of impedance controls for typical tasks encountered in CRM systems.  

Impedance relationships have also been used in conjunction with adaptive control, as an alternative to reduce the system's sensitivity due to modeling mismatch. In \cite{li2017indirect}, an adaptive strategy is used to generate a desired motion trajectory in compliance with desired forces at the end effectors, even when the payload stiffness is unknown. Furthermore, adaptation is used in \cite{ren2014adaptive} to obtain variable damping and stiffness, exhibiting different impedance characteristics according to the forces exerted by the environment.

In this paper, passivity techniques are used to guarantee system stability in the presence of interaction forces and moments between the payload and the robots. 

The remainder of the paper is organized as follows: In Section~\ref{sec_math_pre}, mathematical preliminaries and a brief background on passivity and $\mathcal{L}_2$ stability results are presented. In Section~\ref{sec_modeling}, a model of the CRM is presented. Section~\ref{sec_implement} discusses the semi-active control implementation and associated energy balance equations. The overall control architecture and the impedance control approach are discussed in sections~\ref{Control_Scheme} and \ref{Impedance_control}, respectively. Section~\ref{opt_prob} introduces the optimization and a numerical simulation is presented in section~\ref{simulation}. Finally, conclusions and recommendations for future work are given in section~\ref{dis_future_work}.
\section{Mathematical preliminaries}\label{sec_math_pre}
In this section, we briefly recall some basic properties and classical results of input-output stability and passive systems \cite{van2000l2,vidyasagar2002nonlinear}. Denote $\mathbb{R}_{+}$ as the set of non-negative real numbers, and let $\mathbb{R}^n$ be the $n$-dimensional vector space over $\mathbb{R}$. Define $\Pi$ as the set  of all measurable, real-valued, $n$-dimensional functions of time $f(t):\mathbb{R}_{+}\rightarrow \mathbb{R}^n$ and $t \in \mathbb{R}_{+}$ and define the sets:
\begin{subequations}
\begin{align}
\mathcal{L}^n_2 &\triangleq \{x \in \Pi|~ \| f\|^2_2\triangleq \int_{0}^{\infty} \| f(t)\|^2dt < \infty \}\\
\mathcal{L}^n_{\infty} &\triangleq \{x \in \Pi|~ \| f\|^2_{\infty}\triangleq \sup_t \| f(t)\|^2 < \infty \}
\end{align}
\label{L2}
\end{subequations}
with $\|.\|$ the standard Euclidean norm. Also, extended spaces are defined as:
\begin{subequations}
\begin{align}
\mathcal{L}^n_{2e} &\triangleq \{x \in \Pi|~ \| f\|^2_{2,T}\triangleq \int_{0}^{T} \| f(\xi)\|^2d\xi < \infty,\forall~T \in \mathbb{R}_{+} \}\\
\mathcal{L}^n_{\infty e} &\triangleq \{x \in \Pi|~ \| f\|^2_{\infty,T}\triangleq \sup_{T} \| f(T)\|^2 < \infty,\forall~T \in \mathbb{R}_{+} \}
\end{align}
\label{L2e}
\end{subequations}
Here, $\mathcal{L}^n_{2}$ and $\mathcal{L}^n_{2e}$ are inner product spaces. The truncated inner product of two functions $u(t),y(t) \in \mathcal{L}^n_{2e}$ is defined as:
\begin{equation}
\langle y|u\rangle_T\triangleq \int_{0}^{T} u^T(\xi)y(\xi)d\xi
\label{inner_product}
\end{equation}
where the two functions can be related by a dynamical operator, $H$, such that $y=H(u)$. The operator $H$ is called causal if the value of output, $y$, at time $t$ is a function of the value of input, $u$, up to time $t$ \cite{vidyasagar2002nonlinear}. For a causal operator, we have the following definitions:

\begin{definition} \label{IOSP_def}
Let $H:\mathcal{L}^n_{2e}\rightarrow \mathcal{L}^n_{2e}$ be a causal dynamic operator and $y=H(u)$ is the output. Then for some non-negative constants $\beta$, $\epsilon$ and $\delta$, and $\forall~u\in\mathcal{L}^n_{2e}$, $H(.)$ is called input-output strict passive (IOSP) if:
\begin{equation}
    \langle y|u\rangle_T  \geq  \epsilon\| u\|^2_{2,T} +\delta\| y\|^2_{2,T}-\beta 
    \label{IOSP}
\end{equation}
$H$ is called passive if $\epsilon=\delta=0$, output strictly passive (OSP) if $\epsilon=0$, and input strictly passive (ISP) if $\delta=0$ \cite{vidyasagar2002nonlinear}.
\end{definition}
\begin{definition} \label{L2stability_def}
The causal dynamical operator $H:\mathcal{L}^n_{2e}\rightarrow \mathcal{L}^n_{2e}$ is said to be finite-gain $\mathcal{L}_{2}$ stable if there are non-negative constants $\gamma_0$ and $\beta_0$ such that\cite{ortega2013passivity,vidyasagar2002nonlinear}:
\begin{equation}
    \| y\|_{2,T} \leq \gamma_0\| u\|_{2,T} + \beta_0~~~~~\forall~T\geq 0
\end{equation}

\end{definition}

So far, we have defined the concept of input-output stability for a single causal system. The results can be extended for the feedback system as shown in Fig.~\ref{interconnection}. Here we assume the feedback interconnection is well-defined. That is for $r \triangleq (r_1,r_2)$ and $y \triangleq (y_1,y_2)$, the mapping $y=\tilde{H}(r)$ is causal and $\tilde{H}:\mathcal{L}^n_{2e}\rightarrow \mathcal{L}^n_{2e}$. 

We elaborate the following result for the feedback interconnection: 

\begin{lemma}\label{theorem1}
Consider the well-defined feedback interconnection of two IOSP causal operator dynamical systems $H_1$ and $H_2$ as shown in Fig.~\ref{interconnection}. The following inequality is effective:
\begin{equation} \label{inter_ineq}
c_1\| y_1\|_{2,T}^2+c_2\| y_2\|_{2,T}^2 \leq c_3\| r_1\|_{2,T}^2+c_4\| r_2\|_{2,T}^2-\hat{\beta}
\end{equation}
where:
\begin{align*}
c_1 &= \{(\delta_1+\epsilon_2)-\frac{\gamma_1}{2}-\epsilon_2\gamma_4\}\\
c_2 &= \{(\delta_2+\epsilon_1)-\frac{\gamma_2}{2}-\epsilon_1\gamma_3\}\\
c_3 &=\{\frac{1}{2\gamma_1}+\frac{\epsilon_1}{\gamma_3}\}\\
c_4 &=\{\frac{1}{2\gamma_2}+\frac{\epsilon_2}{\gamma_4}\}
\end{align*}
for $\gamma_j > 0, j\in \{1,2,3,4\}$ and constant $\hat{\beta} >0$.
\begin{figure}
\begin{center}
\includegraphics[width=0.25\textwidth]{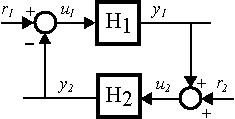}    
\caption{Feedback interconnection with two external inputs.}  
\label{interconnection}                         
\end{center}                                
\end{figure}

\begin{proof}
Since $H_1$ and $H_2$ are IOSP, we have:
\begin{equation}
    \langle y_i|u_i\rangle_T  \geq  \epsilon_i\| u_i\|^2_{2,T} +\delta_i\| y_i\|^2_{2,T}-\beta_i \quad i\in \{1,2\}
    \label{prof_pass_1}
\end{equation}
From Fig.~\ref{interconnection}, we have $u_1 = r_1-y_2$ and $u_2 = r_2+y_1$, therefore:
\begin{equation}
    \langle y_1|u_1\rangle_T+\langle y_2|u_2\rangle_T=\langle y_1|r_1\rangle_T+\langle y_2|r_2\rangle_T
    \label{prof_pass_2}
\end{equation}
Using the interconnection relations, we can write:
\begin{subequations}
\begin{align}
    \| u_1\|^2_{2,T} &= \| r_1\|^2_{2,T}+\| y_2\|^2_{2,T}-2\langle y_2|r_1\rangle_T\notag\\
    &\geq \| y_2\|^2_{2,T}-2\langle y_2|r_1\rangle_T\\
    \| u_2\|^2_{2,T} &= \| r_2\|^2_{2,T}+\| y_1\|^2_{2,T}+2\langle y_1|r_2\rangle_T\notag\\
    &\geq \| y_1\|^2_{2,T}+2\langle y_1|r_2\rangle_T
\end{align}
\label{prof_pass_3}
\end{subequations}
Using equations~(\ref{prof_pass_1}), (\ref{prof_pass_2}), and (\ref{prof_pass_3}), we have the following inequality:
\begin{equation}
\begin{split}
        &\langle y_1|r_1\rangle_T+\langle y_2|r_2\rangle_T+2\epsilon_1\langle y_2|r_1\rangle_T-2\epsilon_2\langle y_1|r_2\rangle_T \geq\\
        &(\delta_1+\epsilon_2)\| y_1\|^2_{2,T}+(\delta_2+\epsilon_1)\| y_2\|^2_{2,T}-\beta_1-\beta_2
\end{split}
    \label{prof_pass_4}
\end{equation}
Using the Cauchy–Schwarz inequality, for any two arbitrary signals $r$ and $y$, there is a $\gamma >0$, such that:

\begin{subequations}
\begin{align}
           \langle y|r\rangle_T &\leq \| y\|_{2,T}\| r\|_{2,T}+\frac{1}{2}(\frac{1}{\sqrt{\gamma}} \| r\|_{2,T}-\sqrt{\gamma}\| y\|_{2,T})^2\notag\\
    &\leq \frac{1}{2\gamma}\| r\|^2_{2,T}+\frac{\gamma}{2}\| y\|^2_{2,T}\\
           \langle y|r\rangle_T &\geq -\| y\|_{2,T}\| r\|_{2,T}-\frac{1}{2}(\frac{1}{\sqrt{\gamma}} \| r\|_{2,T}-\sqrt{\gamma}\| y\|_{2,T})^2\notag\\
    &\geq -\frac{1}{2\gamma}\| r\|^2_{2,T}-\frac{\gamma}{2}\| y\|^2_{2,T}
\end{align}
\label{prof_pass_5}
\end{subequations}
If we assign $\gamma_1,\gamma_2,\gamma_3$ as constants for the pair of signals $(r_1,y_1)$, $(r_2,y_2)$, and $(r_1,y_2)$, respectively, in inequality~(\ref{prof_pass_5}a), and $\gamma_4$ for $(r_2,y_1)$ in inequality~(\ref{prof_pass_5}b), then by substituting the obtained inequalities in inequality~(\ref{prof_pass_4}) and with some simplification and manipulation, equation~(\ref{inter_ineq}) can be obtained. 
\end{proof}
\end{lemma}
Finally, we list two useful standard results:
\begin{lemma}
The feedback in Fig.~\ref{interconnection} is finite-gain $\mathcal{L}_2$ stable if $\delta_1+\epsilon_2 >0$, $\delta_2+\epsilon_1>0$ \cite{brogliato2019dissipative}.
\end{lemma}
\begin{theorem} \label{theorem_2}
Let $y=H(s)u$, where $H(s) \in \mathbb{R}^{n \times n}$ is an exponentially stable and strictly proper transfer function. If $u \in \mathcal{L}_2$, then $y \in \mathcal{L}_2 \cap \mathcal{L}_{\infty}$, $\dot{y} \in \mathcal{L}_2$ and $y \rightarrow 0$ as $t \rightarrow \infty$. In addition, if $u \rightarrow 0$ as $t \rightarrow \infty$, then $\dot{y} \rightarrow 0$ \cite{vidyasagar2002nonlinear,kelly1989adaptive}.
\end{theorem}

\begin{figure}
\begin{center}
\includegraphics[width=0.38\textwidth]{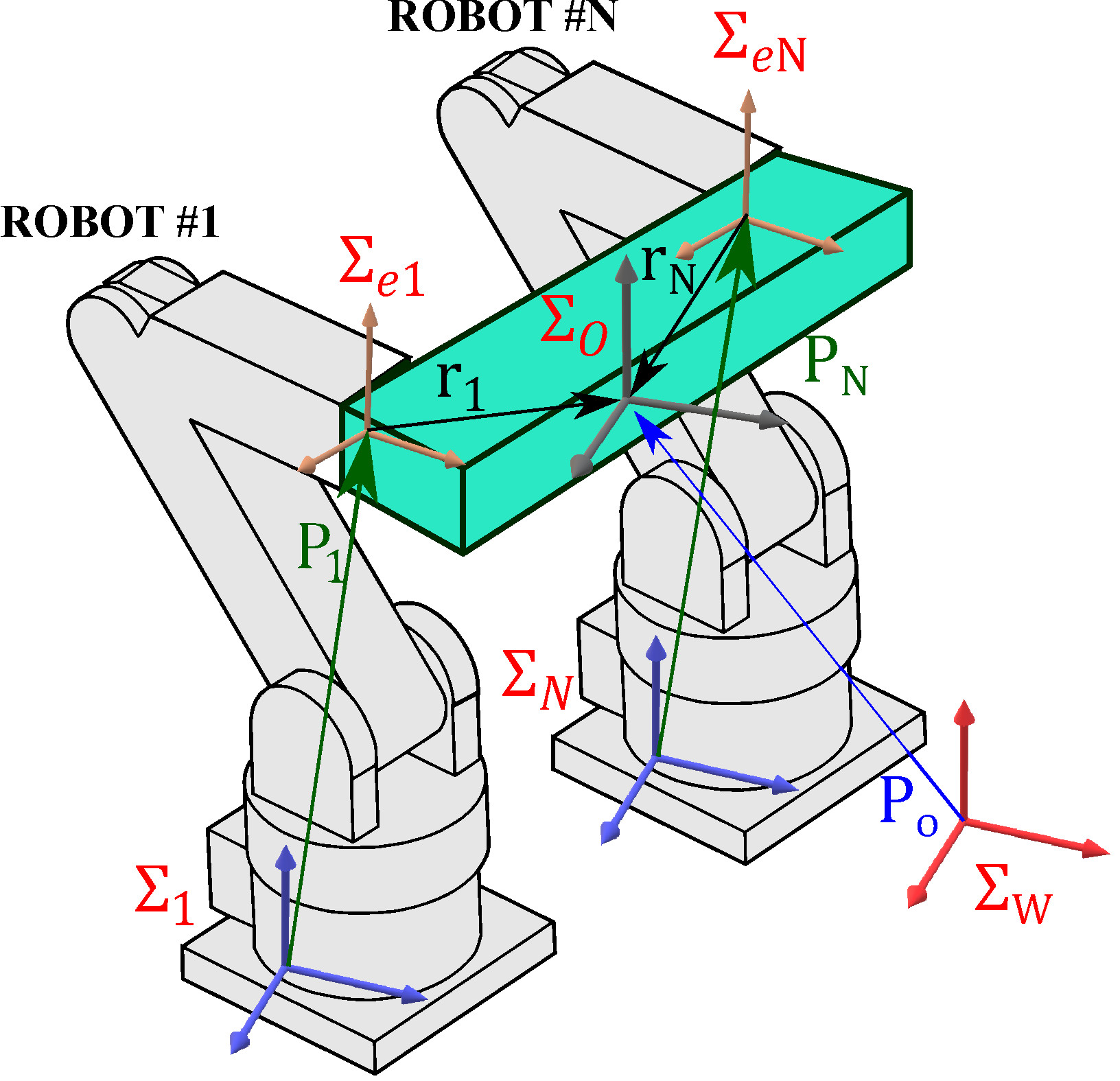}    
\caption{Illustration of $N$ robots grasping an object. The coordinate frames attached to the end-effector, robots' base, load's center of mass, and arbitrary world frame are, respectively, $\Sigma_{ei}$, $\Sigma_{i}$, $\Sigma_{O}$, $\Sigma_{W}$.}  
\label{CRM_RRR}                                 
\end{center}                                
\end{figure}
\section{Modeling}\label{sec_modeling}
In this section, we introduce the dynamical model of $N$ cooperative robots carrying a rigid load by grasping it as shown in Fig.~\ref{CRM_RRR}. The objective in this section is to model the robots, the load and to establish the coupling between the robots. First we build the augmented model of each robot, its joint mechanisms and the DC motor actuators. Then, based on the augmented model of the individual robots, a comprehensive model of the cooperative robots is introduced. The load is modeled assuming to be a rigid object. Finally, the kinematic and dynamic coupling of the closed chain system is discussed.  

To develop the dynamics equations, for the remainder of this paper, unless otherwise told, the following notations are used: any parameter with subscript $i$ indicates the parameter for i\textit{th} robot, $i\in \{1,...,N\}$, and subscript $ij$ indicates the parameter for j\textit{th} joint of i\textit{th} robot, $j\in\{1,2,...,n\}$.
\subsection{Comprehensive dynamics of cooperative robots}
We assume there are $N$ non-redundant robots carrying an object. The dynamics of i\textit{th} robot, without considering the drive systems, modeled as:
\begin{equation}
D^{o}_i(q_i)\ddot{q}_i+C^{o}_i(q_i,\dot{q}_i)\dot{q}_i+g_i(q_i)=\tau_i+J_i^T(q_i)F_i
\label{robot_i_dyn}
\end{equation}
where $q_i~\in~\mathbb{R}^{n}$ is the vector of joint coordinates, $D^{o}_i(q_i)$ is the inertia matrix, $C^{o}_i(q_i)$ is the matrix accounting for centripetal and Coriolis effects, $g_i(q_i)$ is the gravity vector, $J_i~\in~\mathbb{R}^{n \times n}$ is the Jacobian matrix which is assumed to be nonsingular, $F_i=[f^T_i,m^T_i]^T~\in~\mathbb{R}^{n}$ is the vector of forces $(f_i)$ and moments $(m_i)$ applied by the object at the end-effector, and $\tau_i=[\tau_{i1},\tau_{i2},...,\tau_{in}]^T~\in~\mathbb{R}^{n}$ is the vector of forces/moments applied at the joints.

\begin{figure}
\begin{center}
\includegraphics[width=0.33\textwidth]{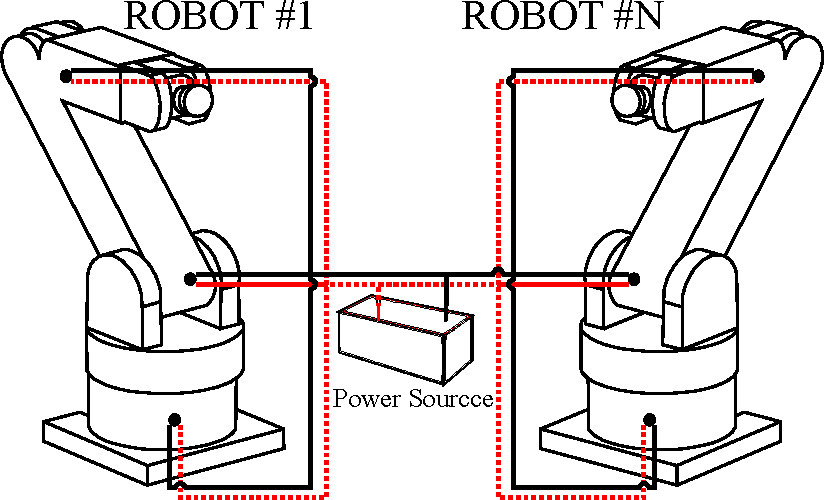}    
\caption{Illustration of semi-active joints connection to a single power source.}  
\label{semiactive_con}                                 
\end{center}                                
\end{figure}

All $n$ joints of the robot are considered to be \emph{semi-active}. A semi-active joint mechanism (JM) only exchange mechanical power with the robot, and it is connected to an electric energy storage element source, e.g. an ultracapacitor. We consider that all semi-active joints of all robots are connected to a common storage element, as illustrated in Fig~\ref{semiactive_con}, and they are regenerative, allowing the power source to be charged whenever surplus mechanical energy from the robots is available to flow back through DC motors. 
Also, in each semi-active JM, the interface torque between the JM and the j\textit{th} robot link can be described by:
\begin{equation}
\tau_{ij}=-J_{ij}\bar{n}^2_{ij}\ddot{q}_{ij}-b_{ij}\bar{n}^2_{ij}\dot{q}_{ij}+\bar{n}_{ij}\tau_{ind_{ij}}
\label{JM}
\end{equation}
where $J_{ij}$, $\bar{n}_{ij}$, and $b_{ij}$ are the JM moment of inertia, the gear ratio, and friction coefficient, respectively. $\tau_{ind_{ij}}=\alpha_{ij}I_{{ij}}$ is the induced torque of the motor, where $I_{{ij}}$ is the current in the motor and $\alpha_{ij}$ is the motor constant.

\begin{figure}
\begin{center}
\includegraphics[width=0.45\textwidth]{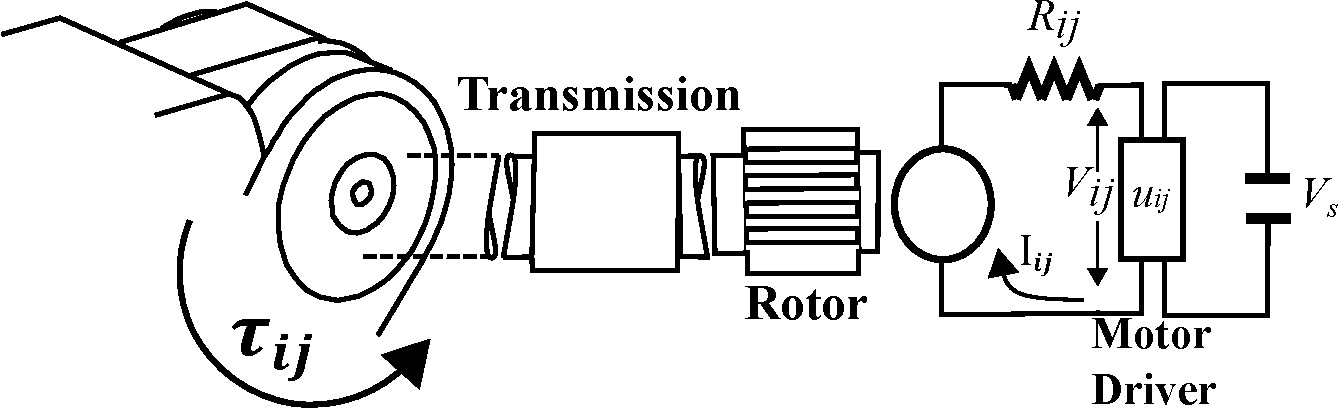}    
\caption{Semi-active joint setup and its connection to motor driver and ultracapacitor.}  
\label{transmission}                                 
\end{center}                                
\end{figure}

By combining the dynamics of robot and JMs, the augmented model of robot-JMs can be obtained. This can be done by finding $\tau_{ind_{ij}}$. In each JM, an ideal regenerative four-quadrant power conversion element (motor driver) is used to control the amount and
direction of the applied voltage to the DC motor. Fig.~\ref{transmission} shows the model of the JM, where the converter voltage ratio is defined as  $u_{ij}=V_{{ij}}/V_{s}$. With this,  the applied voltage to a motor can be written as:
\begin{equation}
V_{{ij}}=R_{ij}I_{{ij}}+a_{ij}\dot{q}_{ij}
\label{Viq}
\end{equation}
where $R_{ij}$ is the resistance of the motor, and $a_{ij}=\alpha_{ij}\bar{n}_{ij}$ \cite{khalaf2019trajectory}.  By substituting $V_{{ij}}$ with $u_{ij}$ in equation~(\ref{Viq}) and obtaining the induced torque of motor as a function of voltage ratio and angular velocity and finally, substituting the induced torque in equation~(\ref{JM}), the applied torque can be rewritten as:
\begin{equation}
\tau_{ij}=-J_{ij}\bar{n}^2_{ij}\ddot{q}_{ij}-(b_{ij}\bar{n}^2_{ij}+\frac{a_{ij}^2}{R_{ij}})\dot{q}_{ij}+\frac{a_{ij}u_{ij}}{R_{ij}}V_{s}
\label{JM_aug}
\end{equation}
Replacing $\tau_{ij}$ from equation~(\ref{JM_aug}) into equation~(\ref{robot_i_dyn}) and moving terms containing $\dot{q}_i$ and $\ddot{q}_i$ from right- to left-hand side, yields the following augmented robot-JM model:
\begin{equation}
D_i(q_i)\ddot{q}_i+C_i(q_i,\dot{q}_i)\dot{q}_i+g_i(q_i)=\mathcal{T}_i+J_i^T(q_i)f_i
\label{robot_i_dynamic_aug}
\end{equation}
where in the coupled dynamics, $D^{o}_i$ and $D_i$ are the same in all elements except in diagonal elements such that: $D_{i}=D^{o}_{i}+\mbox{diag}(J_{ij}\bar{n}^2_{ij})$, also $C_i$ and $C^{o}_i$ are only different in diagonal terms such that $C_{i}=C^{o}_{i}+\mbox{diag}(b_{ij}\bar{n}^2_{ij}+a_{ij}^2/R_{ij})$, and $\mathcal{T}_i=[U_{i1},...,U_{i{n}}]^T$, where:
\begin{equation}
  U_{ij}=\frac{a_{ij}u_{ij}}{R_{ij}}V_{s}
\label{T_i}
\end{equation}
Finally, by combining the augmented dynamics of all $N$ robots, the comprehensive dynamics of CRM can be written as:
\begin{equation}
D(q)\ddot{q}+C(q,\dot{q})\dot{q}+G(q)=\mathcal{T}+J^T(q)F
\label{comp_dyn}
\end{equation}
where $X=[X^T_1,...,X^T_N]^T$, $X \in \{q,G,\mathcal{T},F\}$ and $Y=\mbox{diag}(Y_1,...,Y_N)$, $Y \in \{D,C,J\}$. Hereafter, for ease of notation, we use $J^T(q)F \triangleq \mathcal{T}_{ext}$.

\subsection{Dynamics of the load}
The dynamics and kinematics of manipulators are coupled because the contact forces in all robots are interacting through the load. The dynamic coupling effects can be measured using the dynamics of the load. It is assumed that the mass of load is in its center of mass at $\Sigma_O$ and the orientation is the orientation of coordinate frame $\Sigma_O$ as shown in Fig.~\ref{CRM_RRR}. In the world frame ($\Sigma_W$), the translation and rotation dynamics for the load are:
\begin{equation}
\begin{bmatrix}
           m_oI_{3\times3}&0_{3\times3} \\
           0_{3\times3}&I_o 
\end{bmatrix}
\begin{bmatrix}
           \ddot{P}_o \\
           \dot{\omega}_o
\end{bmatrix}+
\begin{bmatrix}
           -m_o\textbf{g} \\
           \Omega_oI_o\omega_o
\end{bmatrix}=J_o^Tf\triangleq F_o
\label{load_dynamic}
\end{equation}
where $[P_o^T,\omega_o^T]^T=[x_o,y_o,z_o,\omega_x,\omega_y,\omega_z]^T$ is the vector of positions and angular velocities of the load. $I_{3\times3}$ and $0_{3\times3}$ are the identity and the null matrices, respectively, $m_o$ is the mass of the load, and $I_o$ is the moment of inertia about the center of mass. $\Omega_o$ is a skew-symmetric matrix, representing the cross product of $\omega_o\times I_o\omega_o$, $F_o$ is the vector of all forces/moments applied at the center of mass of load, and $J_o^T=[J_{o_1}^T ... J_{o_N}^T]\in\mathbb{R}^{6\times6N}$ is called the grasp matrix and defined as:
\begin{equation}
J_o^T=
\begin{bmatrix}
           I_{3\times3}&0_{3\times3}&...&I_{3\times3}&0_{3\times3} \\
           -S(r_1)&I_{3\times3}&...&-S(r_N)&I_{3\times3}  
\end{bmatrix}
\label{grasp_matrix}
\end{equation}
where $r_i$ is the vector from $\Sigma_{ei}$ to $\Sigma_{o}$ (see Fig.~\ref{CRM_RRR}) and $S(r_i)$ is a skew-symmetric matrix, representing $r_i\times F_i$. The applied forces/moments by manipulators, $f$, can be decomposed to motion-inducing and internal parts as:
\begin{equation}
f=f_M + f_I
\label{force}
\end{equation}
where $f_M$ is the motion-inducing part that may balance the object's dynamics and $f_I$ is the internal part consists of compressive, tensile and torsion forces/moments. Since $f_I$ does not contribute to the motion of the object and has no net force, thus it spans the null space of $J_o^T$. Based on this fact, decomposition as equations~(\ref{force_decom}) can be used to find $f_M$ and $f_I$ where $(J_o^T)^+$ is the generalized inverse of $J_o^T$ \cite{bonitz1996internal}.
\begin{subequations}
\begin{align}
f_M \triangleq &[f_{M_1}^T,...,f_{M_N}^T]^T=(J_o^T)^+ J_o^Tf\\    
f_I \triangleq &[f_{I_1}^T,...,f_{I_N}^T]^T=(I_{6N\times6N}-(J_o^T)^+ J_o^T)f
\end{align}
\label{force_decom}
\end{subequations}
Using equations~(\ref{force_decom}) requires measuring force/moment at the end-effectors. This can be accomplished by installing a force sensor. To avoid direct measurement of force/moment, methods based on load distribution are introduced in \cite{erhart2015internal,walker1991analysis,kawasaki2006decentralized}  
\subsection{Coupling in cooperative robots}
The dynamics (load dynamics) and kinematics coupling cause constraint motion in each robot due to the closed chain configuration of the CRM. The kinematics coupling is defined as a set of equations, relating the translational and the rotational motions of the end-effectors that hold and move the rigid load securely \cite{koivo1990modeling}. In the world frame, the translational and rotational constraints are defined as:
\begin{subequations}
\begin{align}
    P_{i}(q_{i})+r_{i}(q_{i})&=P_{\acute{i}}(q_{\acute{i}})+r_{\acute{i}}(q_{\acute{i}})\\
    \Gamma_i(q_i) - \Gamma_k(q_{\acute{i}}) &= \delta R_{i\acute{i}}\quad i,\acute{i}\in \{1,...,N\}
\end{align}
\label{kin_const}
\end{subequations}
where the first constraint expresses the concatenated vectors, through each robot, that connect $\Sigma_w$ to $\Sigma_o$, $\Gamma_i(q_i)$ is the vector that expresses the orientation of the i\textit{th} end-effector, and $\delta R_{i\acute{i}}$ is a constant vector showing the difference between the orientation of the end-effectors at all time of the maneuver.

For $N$ robots, there are $N(N-1)/2$ translational constraint vectors in the format of equation~(\ref{kin_const}a) which each constraint contains three independent scalar equations, i.e. constraints in $(X,Y,Z)_W$ directions. There are $N(N-1)/2$ rotational constraint vectors in the format of equation~(\ref{kin_const}b). Depending on the representation of the orientation, e.g. Euler or quaternion, various number of scalar constraints can be obtained.
\section{Control implementation and energy balances}\label{sec_implement}
To devise a control scheme for CRM, there are two main consecutive steps. The first step is to design a controller to achieve the motion control objective of CRM; moving an object along a desired trajectory while grasping it securely without any damage to either object or robots. In the second step, a relationship is established between the control signal in the first step and the voltage ratios $u_{ij}$, which are available control inputs. This arises from the incorporation of storage element voltage feedback $V_s$ in the augmented model. The voltage ratios $u_{ij}$ are adjusted in this work using a method called semi-active virtual control (SVC) \cite{richter2015framework}. Using SVC enables us to describe energy exchange with the storage element and write energy balance equations for the closed system including robots and storage elements. 

\subsection{Semi-active virtual control strategy}
Assume that a motion controller, called virtual control law ($\tau^v_{ij}$), for $U_{ij}$ in equation~(\ref{comp_dyn}) has been designed. For $\tau^v_{ij}$ to meet the control objectives for
the augmented model, a solution for $u_{ij}$ is sought that enforces:
\begin{equation}
\tau^v_{ij}=\frac{a_{ij}u_{ij}}{R_{ij}}V_{s}
\label{SVC}
\end{equation}
To design the virtual control, any feedback law compatible
with the desired motion control objectives can be selected. If virtual matching, i.e. $\tau^v_{ij}=U_{ij}$, holds at all
times, properties that apply to the virtual design such as stability, tracking performance, robustness, etc. will be propagated to the actual system. Virtual matching is possible as long as the storage element has nonzero voltage and it will hold exactly whenever the power converters are not saturated (e.g., $-1 \leq u_{ij} \leq 1$) and there is accurate knowledge of parameters $a{ij}$ and $R_{ij}$. The modulation law for exact virtual matching is obtained by solving for $u_{ij}$ from equation~(\ref{SVC}).

The SVC technique uses $V_s$ as feedback information in the virtual matching of equation~(\ref{SVC}). This approach permits the formulation of control laws and energy balance equations \emph{without the need for a model of the storage element}. Devices such as ultracapacitors and batteries have complex and uncertain models. 
\subsection{Internal and external energy balance}
The internal energy balance describes the power exchange
between the semi-active joints and power storage element. Assuming there is no power lost  during power exchange between the power source and DC motors, the input/output energy balance for the motor drives can be written as:
\begin{equation}
V_{s}I_{s}= \sum_{i=1}^{N}\sum_{j=1}^{n} V_{{ij}}I_{{ij}}
\label{power_balance}
\end{equation}
where $V_{s}$ and $I_{s}$ are the voltage and current of the power storage, respectively. Dividing both sides by $V_{s}$, substituting current from equation~(\ref{Viq}) into equation~(\ref{power_balance}), and taking the integral in the time interval $[0,T]$, the internal energy balance equation can be obtained as follows:
\begin{equation}
\Delta\,E_s =\int_{0}^{T}(\dot{q}^T{\mathcal{T}}^{v}-({\mathcal{T}}^{v})^TR_a{\mathcal{T}}^v)dt
\label{internal_e_balance}
\end{equation}
where $\Delta E_s$ is the electric energy change in the power storage, $\mathcal{T}^v=[(\tau^{v}_1)^T,...,(\tau^{v}_N)^T]^T$, $\tau^{v}_i=[\tau^{v}_{i1},...,\tau^{v}_{in}]^T$, and $R_a=\mbox{diag}(R_{ij}/a^2_{ij}) \in \mathbb{R}^{nN\times nN}$. The electric energy, $E_s$ is assumed to be a positive-definite function of the power storage charge $y$, i.e. $E_s=E_s(y)$. We assume during a maneuver, there is enough charge in the storage element for maneuver execution. An example of a power source is the ultracapacitor. All results obtained in this section can be used interchangeably for virtual control design using an ultracapacitor as power source, i.e. $V_s \triangleq V_{cap}$. 

An external energy balance for CRM can be derived as:
\begin{equation}
W_{{ext}} = \Delta\,E_{c}\,+\,\Delta\,E^{Tot}_{m}\,+\,\Sigma^{Tot}_{m}\,+\,\Sigma_{e}
\label{external_balance}
\end{equation}
where $W_{{ext}}$ is the work done by external forces/moments, $\Delta E^{Tot}_m$ is the change in the mechanical energy of CRM, $\Sigma^{Tot}_{m}$ is the mechanical losses, and $\Sigma_{e}$ is the resistance heating (Joule) losses \cite{khalaf2019trajectory}. The derivation and definition of terms in equation~(\ref{external_balance}) is outlined in Appendix~\ref{app_energy}.

\begin{figure*}[h]
\centering
\includegraphics[width=0.8\textwidth]{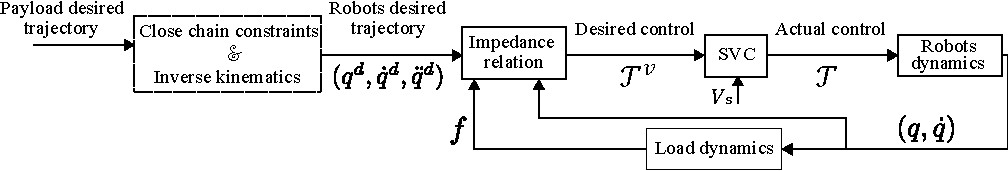}
\caption{\textit{Control scheme to control CRM in joint space, using impedance relation and SVC.}}
\label{control_scheme}
\end{figure*}

\section{Control scheme}\label{Control_Scheme}
In this study, the motion task is defined as following a desired trajectory closely but allowing deviations from the desired trajectory to limit the interaction forces/moments that could damage the robots or load. This can be captured by requiring all internal signals in a feedback control loop to be bounded, which is conveniently achieved with the passivity and $\mathcal{L}_2$ stability tools. Therefore, a control scheme based on impedance relation is developed to accomplish the motion task. Impedance control enforces a compliant behavior to avoid large values of contact force/moment \cite{bonitz1996internal,caccavale2001impedance}. This technique is particularly useful in CRM since there is interaction between robots and the load.

The proposed control scheme is according to Fig~(\ref{control_scheme}). We assume the load's center of mass is following a pre-designed desired trajectory in world space. That is, the translational and rotational information of center of mass are $P^d_o$ and $R^d_o$, respectively.

For the i\textit{th} robot, the end-effector's desired position ($P^d_i$), orientation ($R^d_{ei}$), linear velocity ($\dot{P}^d_i$),  angular velocity ($\omega^d_{ei}$), proportional linear acceleration ($\ddot{P}^d_i$), and angular acceleration ($\dot{\omega}^d_{ei}$) can be calculated using:
\begin{equation}
\begin{split}
P^d_i&=P^d_o -R^d_or_i \\
R^d_{ei}&=R^d_o  \\
\dot{P}^d_i&=\dot{P}^d_o +S(R^d_or_i)\omega^d_o \\
\omega^d_{ei}&=\omega^d_o \\
\ddot{P}^d_i&=\ddot{P}^d_o +S(\omega^d_o)S(R^d_or_i)\omega^d_o+S(R^d_or_i)\dot{\omega}^d_o \\
\dot{\omega}^d_{ei}&=\dot{\omega}^d_o \\
\end{split}
\label{end_effector_des}
\end{equation}
where $S(.)$ is the matrix performing cross-product and $R^d_o$ is a matrix transformation that transfers $r_i$ from object's frame ($\Sigma_o$) to the frame attached to the i\textit{th} end-effector ($\Sigma_{ei}$) \cite{caccavale2001impedance}. Using translational and rotational information of the end-effector, from Equation.~(\ref{end_effector_des}), kinematics of the robot, and any inverse kinematics technique, the desired position ($q^d_i$), velocity ($\dot{q}^d_i$) and acceleration ($\ddot{q}^d_i$) of joints, can be obtained \cite{spong2006robot,caccavale2001impedance}. The desired states in joint space along with the actual states and the interaction force between robots and the load are used in the impedance relation to calculate the virtual control, ${\mathcal{T}}^{v}$. The SVC uses the power source voltage information to regulate the virtual control and creates the actual control signal.

\section{Impedance control}\label{Impedance_control}
The characteristics of an impedance control problem are given as a relationship between a desired trajectory and a desired dynamics. In this sense, the system is forced to establish a mathematical relationship between the interaction forces and the position error. To design the virtual control, we enforce the following relationship introduced in \cite{kelly1989adaptive}:

\begin{equation}
\lim_{t\to\infty} M\ddot{\tilde{q}}+B\dot{\tilde{q}}+K\tilde{q}=\mathcal{T}_{ext}
\label{desired_impedance}
\end{equation}
where $\tilde{q}=q^d-q$ and the inertia ($M$), damping ($K$), and stiffness ($K$) are diagonal positive-definite matrices called gain matrices here, and all parameters are defined for each robot accordingly.

To achieve the objective in equation~\ref{desired_impedance}, we define an auxiliary error as:
\begin{equation}
\zeta = \tilde{q}-[p^2M+pB+K]^{-1}\mathcal{T}_{ext}
\label{aux_error}
\end{equation}
where $p=d/dt$. It is easy to verify that if $\zeta$ converges to zero, the desired impedance in equation~(\ref{desired_impedance}) is achieved. Therefore, we define the following variables:
\begin{subequations}
\begin{align}
&S = -(\dot{\zeta}+\Lambda \zeta)\\    
&\dot{q}_r = \dot{q}-S
\end{align}
\label{s_qr}
\end{subequations}
where $\Lambda$ is a diagonal positive definite matrix. The virtual torque can be calculated using:
\begin{equation}
{\mathcal{T}}^{v} = D(q)\ddot{q}_r+C(q,\dot{q})\dot{q}_r+G(q)-K_DS-\mathcal{T}_{ext}
\label{des_tau}
\end{equation}
where $K_D$ is a positive definite matrix. From a practical point of view, $\dot{q}_r$ and $\ddot{q}_r$ in equation.~(\ref{des_tau}) can be implemented by:
\begin{subequations}
\begin{align}
\dot{q}_r &= \dot{q}^d+\Lambda {\tilde{q}} - (pI+\Lambda)[p^2M+pB+K]^{-1}\mathcal{T}_{ext}\\
\ddot{q}_r &= \ddot{q}^d + \Lambda \dot{\tilde{q}} - p(pI+\Lambda)[p^2M+pB+K]^{-1}\mathcal{T}_{ext}
\end{align}
\label{qr_qdotr}
\end{subequations}
which shows that the implementation of the controller does not require the measurement of the acceleration. 

\begin{theorem} \label{prop_main}
Consider the auxiliary error in equation~(\ref{aux_error}) and $S$ in equation~(\ref{s_qr}a). Then the followings hold:
\begin{enumerate}[label=\textbf{(\alph*)}]
    \item $S \in \mathcal{L}_2$
    \item $\zeta \in \mathcal{L}_2 \cap \mathcal{L}_{\infty}$, $\dot{\zeta} \in \mathcal{L}_2$
    \item  $\zeta \rightarrow 0$ as $t \rightarrow \infty$.
    \item $q,\Dot{q} \in \mathcal{L}_{2e}$ and bounded.
\end{enumerate}
and therefore, the closed-loop system is $\mathcal{L}_2$ finite gain stable.  
\begin{proof}
\textbf{(a)}~Assuming ${\mathcal{T}}^{v}={\mathcal{T}}$, by substituting equation~(\ref{des_tau}) in equation~(\ref{comp_dyn}), the closed-loop system is obtained as:
\begin{equation}
D(q)\dot{S}+[C(q,\dot{q})+K_D]S=0
\label{closed_loop}
\end{equation}
We choose the following Lyapunov function candidate for the closed-loop system:
\begin{equation}
V(t) = \frac{1}{2}S^TD(q)S
\label{closed_lyap}
\end{equation}
Differentiating Lyapunov function with respect to time and with some manipulation and simplification, it yields:
\begin{equation}
\dot{V}(t) = -S^TK_DS
\label{lyap_diff}
\end{equation}
Equations~(\ref{closed_lyap}) and(\ref{lyap_diff}) imply that $S$ is asymptotically stable and therefore $S \in \mathcal{L}_2$. 

\noindent \textbf{(b,c)}~Using the definition of $S$ in equation~(\ref{s_qr}a), $\zeta$ and its derivative can be written as:
\begin{subequations}
\begin{align}
\zeta&= -(pI+\Lambda)^{-1}S\\
\dot{\zeta}&= -p(pI+\Lambda)^{-1}S
\end{align}
\label{prop1}
\end{subequations}
where $(pI+\Lambda)^{-1}$ describes a strictly proper, exponentially stable transfer function \cite{ortega2013passivity}. Using Theorem~\ref{theorem_2} and the fact that $S \in \mathcal{L}_2$, we conclude $\zeta \in \mathcal{L}_2 \cap \mathcal{L}_{\infty}$, $\dot{\zeta} \in \mathcal{L}_2$, and hence $\zeta \rightarrow 0$ as $t \rightarrow \infty$. The fact that auxiliary error, $\zeta$, converges to zero, shows the objective of achieving the desired impedance relation in equation~(\ref{desired_impedance}).

\noindent \textbf{(d)}~By taking the derivative of both sides of equation~(\ref{aux_error}) and rearranging the terms, we have:
\begin{equation}
\dot{q}=\dot{q}^d - \dot{\zeta} -p[p^2M+pB+K]^{-1}\mathcal{T}_{ext}
\label{prop2}
\end{equation}
Using equation~(\ref{prop1}b), equation~(\ref{prop2}) can be represented in a feedback interconnection form as depicted in Fig~\ref{feed_trans}(a). Here we define $\mbox{ENV}(\dot{q})$ as a passive operator that maps $\dot{q}$ to $\mathcal{T}_{ext}$. The environment (load) is passive because the load is assumed to be a pure mass that is held by grasping. Since it is passive, we can write:
\begin{equation}
  \langle \mathcal{T}_{ext}|\dot{q}\rangle_T  \geq -\Tilde{\beta}
    \label{ENV}
\end{equation}

\begin{figure}
\begin{center}
\includegraphics[width=0.45\textwidth]{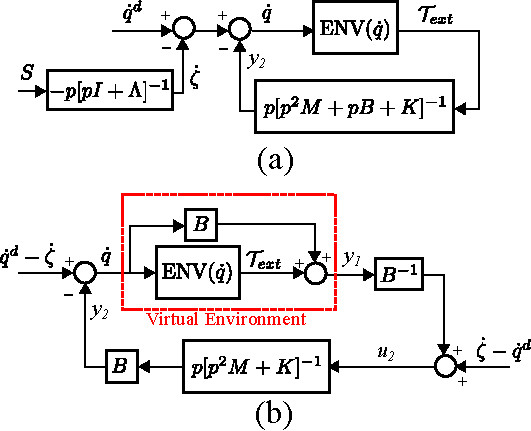}    
\caption{(a) Feedback interconnection with one external input. (b) Transformation of interconnection to two external inputs.}  
\label{feed_trans}                                 
\end{center}                                
\end{figure}
It can be easily shown that the operator $\bar{H}_2 \triangleq p[p^2M+pB+K]^{-1}$ is only passive. Therefore, this feedback loop does not reveal any information about boundedness of signals and it is inconclusive. By performing the loop transformation technique \cite{vidyasagar2002nonlinear,van2000l2}, the equivalent interconnection can be presented as Fig~\ref{feed_trans}(b). In this figure, the virtual environment operator ($H_1$) is composed of passive environment as represented by equation~(\ref{ENV}) and a feed forward term. So it is an ISP system with the following property:
\begin{equation}
  \langle y_1|\dot{q}\rangle_T  \geq -\Tilde{\beta} + \underline{\lambda}_B\| \dot{q}\|^2_{2,T}
    \label{virt-ENV}
\end{equation}
where $\underline{\lambda}_B>0$ is the smallest eigenvalue of the gain matrix $B$. The operator $H_2 \triangleq pB[p^2M+K]^{-1}$, which maps $u_2$ to $y_2$, is only passive. By comparing Fig~\ref{feed_trans}(b) and Fig~\ref{interconnection}, and using equation~(\ref{IOSP}), we can conclude that $\epsilon_1 = \underline{\lambda}_B$, $\epsilon_2=\delta_1=\delta_2=0$. Based on Lemma~\ref{theorem1}, we can write:
\begin{equation}
\bar{c}_2\| y_2\|_{2,T}^2 \leq \bar{c}_3\| \dot{q}^d-\dot{\zeta}\|_{2,T}^2+\bar{c}_4\| \dot{\zeta}-\dot{q}^d\|_{2,T}^2-\bar{\bar{\beta}} 
    \label{prop3}
\end{equation}
where $\bar{c}_i$ can be always chosen as positive numbers. It was shown that $\zeta \in \mathcal{L}_2$ and also using the fact that desired trajectory is bounded, we can conclude $\dot{q}^d,\dot{\zeta}\in \mathcal{L}^n_{2e}$. Therefore, the right-hand side of equation~(\ref{prop3}) is in $\mathcal{L}^n_{2e}$ and so is the left-hand side, i.e. $y_2 \in \mathcal{L}^n_{2e}$ and bounded. Consequently, all signals inside the loop are bounded, i.e. $\dot{q},\mathcal{T}_{ext},y_1,u_2 \in \mathcal{L}^n_{2e}$.

In particular, $\exists~\bar{\mathcal{T}}$ such that $||\mathcal{T}_{ext}|| \leq \bar{\mathcal{T}}$. The impedance relation as a function of the auxiliary error in equation~(\ref{aux_error}), can be rewritten as:
\begin{equation}
M\ddot{Z}+B\dot{Z}+KZ=\mathcal{T}_{ext}
\label{prop4}
\end{equation}
where $Z=\tilde{q}-\zeta$. It can be shown that all solutions of equation~(\ref{prop4}) comply with the inequality:
\begin{equation}
  \mbox{max}(\| \dot{Z}\|_{2,T},\| Z\|_{2,T}) \leq 
  C^{-\kappa t}\sqrt{\| \dot{Z}_0\|^2+\| Z_0\|^2}+\frac{C\bar{\mathcal{T}}}{\kappa}
\label{prop5}
\end{equation}
where $C$ is a constant, $\kappa>0$, and $\dot{Z}_0,Z_0$ are the initial conditions \cite{vukobratovic1997control}. Because the right-hand side of equation~(\ref{prop5}) is monotonically decreasing and tends to a finite value and $\zeta \in \mathcal{L}_2$, we have $q \in \mathcal{L}^n_{2e}$.

\end{proof}
\end{theorem} 

\section{Optimization problem}\label{opt_prob}
We aim to find the impedance gains, $B$ and $K$, that maximize the energy regeneration given the pre-designed references used in the impedance relation of~\ref{desired_impedance}. Maximizing $\Delta E_s$ in equation~(\ref{internal_e_balance}) means less amount of energy is drawn from the storage element. A value of $\Delta E_s>0$ indicates energy regeneration and $\Delta E_s<0$ indicates energy consumption. Considering constant impedance  parameters first, we define the damping and stiffness gains as:
\begin{equation}
\begin{cases}
B = B_c+\bar{B} \\
K = K_c+\bar{K} \\
\end{cases}
\label{gains_cte}
\end{equation}
where $B_c$ and $K_c$ are fixed and $\bar{B}$ and $\bar{K}$ are limited within an upper- and lower-bounds. The optimization is aimed to find $\bar{B}$ and $\bar{K}$ such that $\Delta E_s$ is maximized. Considering equation~(\ref{internal_e_balance}) and the constraints in equation~(\ref{kin_const}), a static optimization problem is formulated as:
\begin{subequations}
\begin{align}
&\max_{\acute{B},\acute{K}}~\Delta E_s=\int_{0}^{T}(\dot{q}^T{\mathcal{T}}^{v}-({\mathcal{T}}^{v})^TR_a{\mathcal{T}}^v)dt\\
&\mbox{subject to:}
\begin{cases}
\mbox{I}&: P_{i}(q_{i})+r_{i}(q_{i})=P_{\acute{i}}(q_{\acute{i}})+r_{\acute{i}}(q_{\acute{i}})\\
\mbox{II}&: \Gamma_i(q_i) - \Gamma_k(q_{\acute{i}}) = \delta R_{i\acute{i}}\quad i,\acute{i}\in \{1,...,N\}
~\\
\mbox{III}&: -V_{s}\bar{a}_R \leq \mathcal{T}^v \leq V_{s}\bar{a}_R\\
\mbox{IV}&: LB_B \leq \acute{B} \leq UB_B\\
\mbox{V}&: LB_K \leq \acute{K} \leq UB_K\\
\end{cases}
\end{align}
\label{opt_form}
\end{subequations}
where the conditions $\mbox{I}$ and $\mbox{II}$ are the translational and rotational constraints, condition $\mbox{III}$ is derived from equation~(\ref{SVC}) and the fact that $-1 \leq u_{ij} \leq 1$. $\bar{a}_R=\mbox{Vec}(a_{ij}/R_{ij}) \in \mathbb{R}^{nN}$ and $\mbox{Vec}(.)$ represents the vectorization operator which creates a
column vector by stacking $a_{ij}/R_{ij}$ for all joints. The lower- and upper-bound for the damping and stiffness gains are indicated in conditions $\mbox{IV}$ and $\mbox{V}$, respectively. These bounds determine the search area for the gains. The problem can be solved with a variety of suitable methods, such as the genetic algorithm used in the simulation example.
\section{Simulation}\label{simulation}
\begin{figure}[h]
\centering
\includegraphics[width=0.35\textwidth]{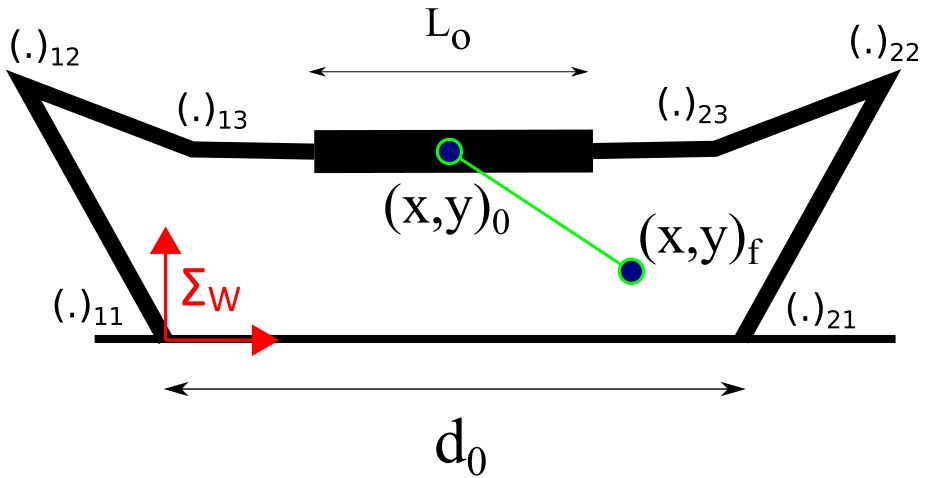}
\caption{\textit{Setup for two identical cooperative robots carrying a load. The world frame is attached to the base of the left robot. The distance between robots is $d_0=0.8~m$ and robots are placed on the same level. $(.)_{ij}$ represents any parameter for the j\textit{th} joint of i\textit{th} robot.}}
\label{RRR_cop}
\end{figure}

To simulate CRM system using the proposed impedance control, two identical RRR planar robots are considered. The robots are grasping a rigid rod as shown in Fig.~\ref{RRR_cop}. It is assumed the actuators in the robots are DC motors and all joints are powered using an ultracapacitor as storage element. The robots' and the rod's parameters are given in Table~\ref{RRR_param}. 

\begin{table}[b]
\small
\begin{center}
\caption{The parameters for two identical robots and the load (rod).}\label{RRR_param}
\renewcommand{\arraystretch}{1.05}
\setlength{\arrayrulewidth}{1.5pt}
\begin{tabular}{lcc}
\hline 
\textbf{Parameter} & \textbf{Value} & \textbf{Unit} \\ 
\hline
\textbf{\underline{Robot:}}&&\\
~~1st, 2nd, and 3rd arm lengths &[0.425,0.39,0.13]&m\\ 
~~1st, 2nd, and 3rd arm masses &[8.05,2.84,1.37]&kg\\ 
~~DC motor resistance &0.4&$\Omega$\\ 
~~DC motor torque constant &0.07&$\Omega$\\
~~DC motor gear ratio &50&-\\ 
\textbf{\underline{Rod:}}&&\\
~~~Mass &5&kg\\
~~~Length ($L_0$)&0.5&m\\
\hline 
\end{tabular}
\renewcommand{\arraystretch}{1}
\end{center}
\end{table} 

The motion task objective is defined as moving the load's center of mass (CM) from initial position, $(x,y)_0$, to final position, $(x,y)_f$, along a pre-designed desired trajectory. The final condition is relaxed to allow for a  search of the optimal impedance gains. Accordingly, a set has been defined as:
\begin{equation}
IS = \{p_o(t)|~||p_o(t)-p^d_o(t_f)||\leq \epsilon_f\}
\label{invar_set}
\end{equation}
where $p_o(t)$ is the position of the load's CM, $p^d_o(t_f)$ is the final desired position, $||.||$ is the Euclidean norm, and $\epsilon_f$ is a scalar boundary. It is assumed that the motion tasks starts from the initial position, and it is accomplished whenever the final position enters the $IS$ set. After finishing the maneuver, the load will be detached from the robots.

The desired load trajectory is based on quintic
polynomial and consists of moving the load 0.4 m in the
x direction, -0.4 m in the y direction, and a rotation of zero degrees in 1 s. The desired trajectories in joint space is obtained using equations~\ref{end_effector_des}.

The fixed inertia, damping, and stiffness impedance gains are selected as $M=\mbox{diag}(18)$, $B_c=\mbox{diag}(197.5)$, and $K_c=\mbox{diag}(825)$. Note that the same parameters are selected for both robots. The optimization in equation~\ref{opt_form} is solved using genetic algorithm (GA). The optimization searches for 12 variables (one damping and one stiffness gain for each joint). The upper- and lower-bounds for the damping and stiffness gains in each joint are $-22\leq \bar{B}_{ij} \leq 22$ and $-75\leq \bar{K}_{ij} \leq 75$, respectively. The values and types of GA operators of this study are given in Table~\ref{ga_para}.

\begin{table}[t]
\small
\begin{center}
\caption{Parameters of genetic algorithm.}\label{ga_para}
\renewcommand{\arraystretch}{1.1}
\begin{tabular}{lc}
\hline 
Parameter & Value  \\ \hline 
Initial population size &50\\ 
Max. No. of generation&30\\ 
Crossover probability&0.75\\ 
Mutation rate&0.02\\ \hline 
\end{tabular}
\renewcommand{\arraystretch}{1}
\end{center}
\end{table} 

\begin{figure}[h]
\centering
\includegraphics[width=0.4\textwidth]{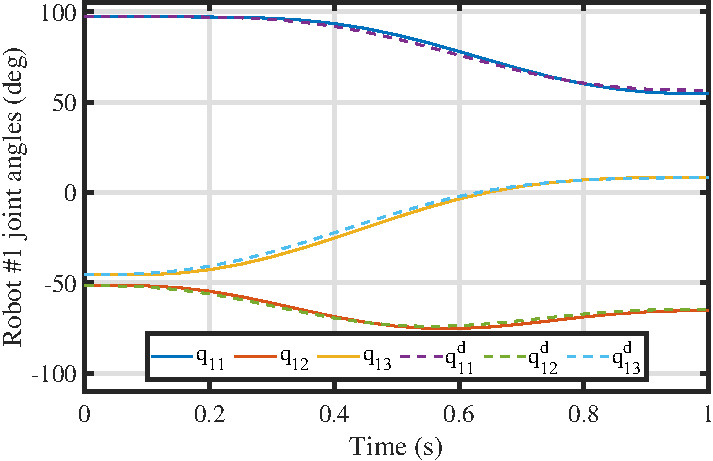}
\caption{\textit{The actual and desired joint angles for Robot \#1.}}
\label{jointR1}
\end{figure}

\begin{figure}[h]
\centering
\includegraphics[width=0.4\textwidth]{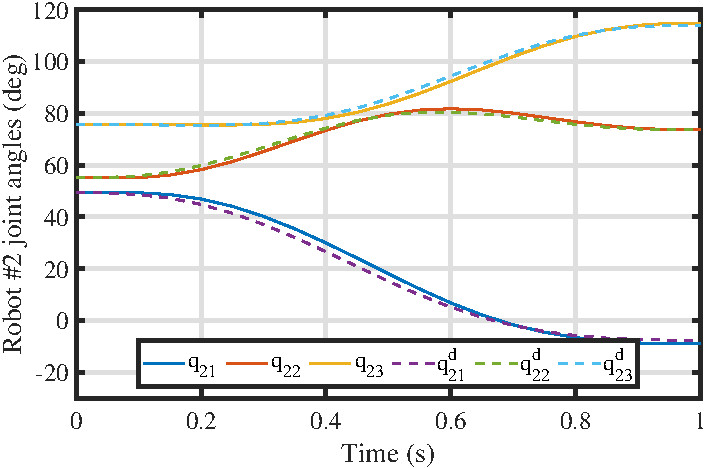}
\caption{\textit{The actual and desired joint angles for Robot \#2.}}
\label{jointR2}
\end{figure}

\begin{figure}[h]
\centering
\includegraphics[width=0.47\textwidth]{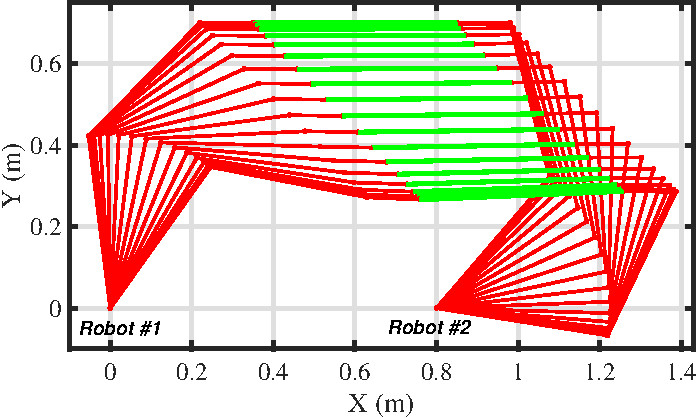}
\caption{\textit{Movement in 2D. Frames are captured in each 0.05 s.}}
\label{plane_mov}
\end{figure}

\begin{figure}[h]
\centering
\includegraphics[width=0.47\textwidth]{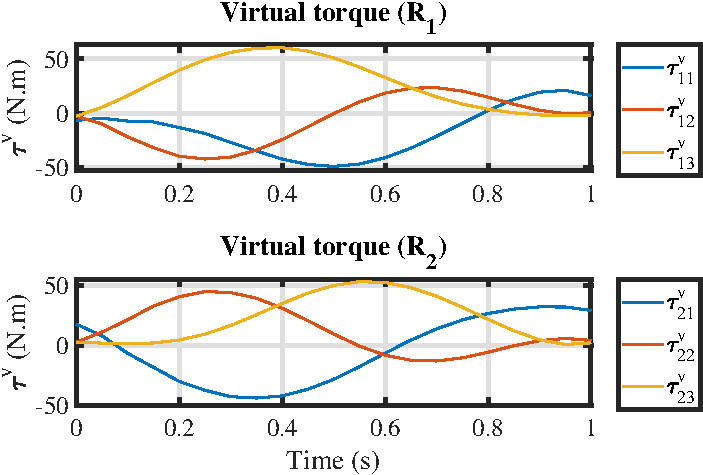}
\caption{\textit{Virtual torques applied to the robots.}}
\label{tau_R1_R2}
\end{figure}

The optimization gives the following impedance gains as optimal values:
\begin{align*}
B_1 = 
    \begin{bmatrix}
     176.1&0&0  \\
     0&178.0&0 \\
     0&0&181.0 
\end{bmatrix}
,B_2 = 
    \begin{bmatrix}
     175.6&0&0  \\
     0&176.8&0 \\
     0&0&176.3
\end{bmatrix}\\
K_1 = 
    \begin{bmatrix}
     752.9&0&0  \\
     0&754.2&0 \\
     0&0&763.6 
\end{bmatrix}
,K_2 = 
    \begin{bmatrix}
     756.4&0&0  \\
     0&761.2&0 \\
     0&0&755.1 
\end{bmatrix}
\end{align*}
where subscripts 1 and 2 denote the matrix for the first and second robot (in Fig.~\ref{RRR_cop}, the robot on the left is named robot 1 and the right one is  called robot 2).
Figures~\ref{jointR1} and~\ref{jointR2} show the time histories of the reference trajectory and the actual angles for both robots. The reference trajectories are closely tracked due to the relatively high value chosen for $K_c$. The 2D movement is depicted in Fig.~\ref{plane_mov} and the virtual controls are shown in Fig.\ref{tau_R1_R2}. 

\begin{figure}[h]
\centering
\includegraphics[width=0.47\textwidth]{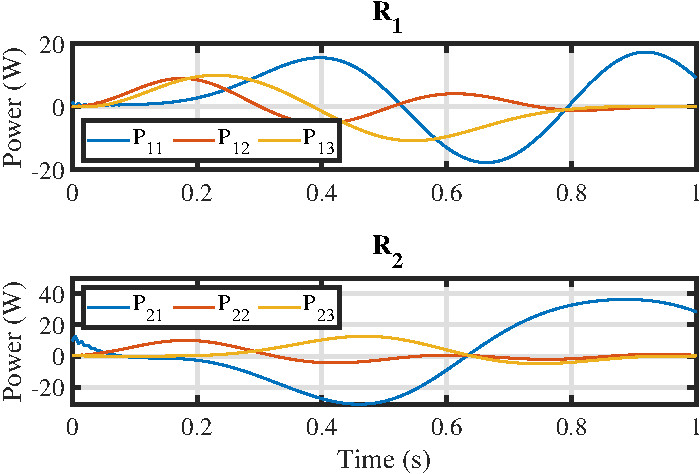}
\caption{\textit{Power in each joint.}}
\label{power_exchange}
\end{figure}

The power consumption in each joint is shown in Fig.~\ref{power_exchange}. Positive power indicates power consumption by the joint and negative power shows energy regeneration. The energy consumption and regeneration is more pronounced in the first joint of each robot.

\begin{figure}[h]
\centering
\includegraphics[width=0.47\textwidth]{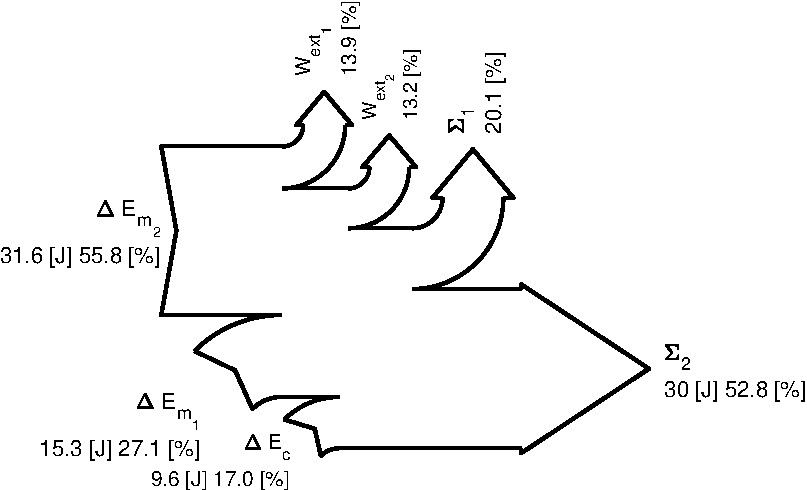}
\caption{\textit{Sankey diagram showing the energy balance for the CRM when following the desired trajectory.}}
\label{sankey_diag}
\end{figure}

Figure~\ref{sankey_diag} shows the Sankey diagrams for the external energy balance of equation~\ref{external_balance}. Most of the energy needed to accomplish the motion task was recovered from the potential energy difference between the initial and final positions. This is a direct consequence of maximizing energy regeneration. To study the effect of energy regeneration, we conclude the results by defining the effectiveness of energy regeneration as \cite{khalaf2019trajectory}:
\begin{equation}
\epsilon = 1-\frac{\Delta E_R}{\Delta E_{NR}}
\label{effective}
\end{equation}
where $\Delta E_R$ and $\Delta E_{NR}$ are the system energy consumption with and without energy regeneration, respectively. $\Delta E_{NR}$ is computed by integrating the power flows in all joints, assuming any negative
power is dissipated (i.e.~$P_{ij}~(P_{ij}~\leq~0)=0$). We have $0\leq \epsilon \leq 1$ where $\epsilon=0$ means energy regeneration
has zero effect in reducing the energy consumption and $\epsilon=1$ indicates that energy regeneration completely reduces energy consumption. For the condition of simulation, $\Delta E_R=9.69~J$ and $\Delta E_{NR}=25.84~J$ results in $\epsilon=0.62$. This shows approximately 60\% reduction in energy consumption due to energy regeneration.

\subsection{Variable impedance parameters}
The optimization in equation~\ref{opt_form} searched for constant values of impedance gain during the motion task. It is possible to look for variable gains, i.e. $B(t)$ and $K(t)$. Variable impedance shows better versatility when there is interaction with environment, and it has been shown to have better performance \cite{ferraguti2013tank,he2016dual,zheng2018time,kishi2003passive,kronander2016stability}. Optimal variable gains may result in further reductions in energy consumption.

If the damping and stiffness gains are time-varying matrices constrained to be positive definite, the relation between $y_2$ and $u_2$ in Fig.~\ref{feed_trans}(b) can be written in the format of a linear time-varying system as:
\begin{subequations}
\begin{align}
    \dot{\mathcal{X}}(t) &=\mathfrak{A}(t)\mathcal{X}(t)+\mathfrak{B}u_2(t)\\
     y_2(t) &= \mathfrak{C}(t)\mathcal{X}(t)
\end{align}
\label{linear}
\end{subequations}
where
$$\mathfrak{A}(t)=\begin{bmatrix}
    \textbf{0} & \textbf{1}\\
    -M^{-1}K(t) & \textbf{0}
    \end{bmatrix},\mathfrak{B}=
    \begin{bmatrix}
    \textbf{0} \\
    M^{-1}
    \end{bmatrix},\mathfrak{C}(t)=
        \begin{bmatrix}
    \textbf{0} & B(t)\\
    \end{bmatrix}$$
In the above matrices, $\textbf{0} \in \mathbb{R}^{nN\times nN}$ is a null matrix and $\textbf{1}=\mbox{diag}(1) \in \mathbb{R}^{nN\times nN}$ is a diagonal matrix.

Variable gains may result in the loss of passivity in the desired impedance relation \cite{ferraguti2013tank}. However, we establish a result providing conditions to preserve passivity in the time-varying case:

\begin{theorem}
The states of the system are bounded and passivity is guaranteed if the following condition is satisfied:
\begin{equation}
\overline{\lambda} < 2\underline{\lambda}^2_B
    \label{pass_const}
\end{equation}
where $\underline{\lambda}_B$ is the smallest eigenvalue of $B(t)$, $\overline{\lambda}$ is the largest eigenvalue of $\dot{P}(t)+P(t)\mathfrak{A}(t)+\mathfrak{A}^T(t)P(t)$, and $P(t)=P^T(t)$ is a positive-definite matrix defined as:
\begin{equation}
 P(t)=
    \begin{bmatrix}
    K(t) & \textbf{0}\\
    \textbf{0} & MB(t)
    \end{bmatrix}
    \label{P-definition}
\end{equation}

\begin{proof}
We assume $B(t)$ and $K(t)$ are bounded, diagonal, positive definite matrices, therefore $P(t)$ is a continuous and bounded positive-definite matrix. Henceforth, to be concise, we will neglect writing the temporal argument. To study the passivity of the time varying linear system in equation~(\ref{linear}),  the following storage function and its temporal derivative can be considered \cite{forbes2010passive}:
\begin{subequations}
\begin{align}
    \mathcal{V}&=\frac{1}{2}\mathcal{X}^TP\mathcal{X}\\
    \dot{\mathcal{V}}&=\frac{1}{2}\mathcal{X}^T\{\dot{P}+P\mathfrak{A}+\mathfrak{A}^TP\}\mathcal{X}+\mathcal{X}^TP\mathfrak{B}u_2
\end{align}
\label{lyap2}
\end{subequations}
Integrating $\dot{\mathcal{V}}$ in the time interval $[0,T]$ gives:
\begin{equation}
\int_{0}^{T}\dot{\mathcal{V}}(\xi)d\xi=\mathcal{V}(t)-\mathcal{V}(0) \geq -\mathcal{V}(0)
    \label{Vdot_int}
\end{equation}
Substituting equation~(\ref{lyap2}b) in the above inequality results in the following inequality:
\begin{equation}
    \int_{0}^{T}[\frac{1}{2}\mathcal{X}^T\{\dot{P}+P\mathfrak{A}+\mathfrak{A}^TP\}\mathcal{X}+\mathcal{X}^TP\mathfrak{B}u_2]d\xi \geq -\mathcal{V}(0)
\end{equation}
Because $P\mathfrak{B}=\mathfrak{C}^T$, assuming $\bar{P}=\dot{P}+P\mathfrak{A}+\mathfrak{A}^TP$, we have:
\begin{equation}
    \int_{0}^{T}\frac{1}{2}\mathcal{X}^T\bar{P}\mathcal{X}d\xi +
    \int_{0}^{T}y^T_2(\xi)u_2(\xi)d\xi \geq -\mathcal{V}(0)
    \label{ineq_v}
\end{equation}
where:
\begin{equation}
    \bar{P}=
     \begin{bmatrix}
    \dot{K} &K-K^TB^T\\
    K^T-BK & M\dot{B}
    \end{bmatrix}
    \label{K_dot}
\end{equation}
To study the passivity of input/output in the above inequality, we study three cases of $\bar{P}$ being negative-definite (ND), positive-definite (PD),  and indefinite (ID).

\noindent 1) In case of ND, the integral including $\bar{P}$ is negative, so inequality~(\ref{ineq_v}) can be written as:
\begin{equation}
  \langle y_2|u_2\rangle_T  \geq -\mathcal{V}(0)
    \label{P_ND_pass}
\end{equation}
which proves the passivity of the system and therefore, Theorem~\ref{prop_main} and its results are valid.

\noindent 2) In the case of PD, because $\bar{P}>0$, passivity can not be concluded directly from equation~(\ref{ineq_v}). In this case,  we consider the following inequalities:
\begin{equation}
    \frac{\underline{\lambda}}{2}\int_{0}^{T}\mathcal{X}^T\mathcal{X}d\xi \leq \int_{0}^{T}\frac{1}{2} \mathcal{X}^T\bar{P}\mathcal{X}d\xi \leq \frac{\overline{\lambda}}{2}\int_{0}^{T}\mathcal{X}^T\mathcal{X}d\xi
    \label{ineq_P_PD}
\end{equation}
where $ \underline{\lambda}>0$ and $\overline{\lambda}$ are the smallest and biggest eigenvalues of $\bar{P}$, respectively.  Also, in the view of output in equation~(\ref{linear}b), we can write:
\begin{equation}
 \underline{\lambda}_B\int_{0}^{T}\mathcal{X}^T\mathcal{X}d\xi  \leq \int_{0}^{T} \mathcal{X}^T\mathfrak{C}^T\mathfrak{C}\mathcal{X}d\xi \leq \overline{\lambda}_B\int_{0}^{T}\mathcal{X}^T\mathcal{X}d\xi
    \label{ineq_P_PD2}
\end{equation}
where $ \underline{\lambda}_B>0$ and $\overline{\lambda}_B$ are the smallest and biggest eigenvalues of $\mathfrak{C}$, respectively. From equations~(\ref{ineq_P_PD}) and (\ref{ineq_P_PD2}), the following inequality can be resulted:
\begin{equation}
 \frac{\overline{\lambda}}{2}\int_{0}^{T}\mathcal{X}^T\mathcal{X}d\xi \leq \frac{\overline{\lambda}}{2\underline{\lambda}_B}\int_{0}^{T} y^T_2y_2d\xi
    \label{ineq_P_PD3}
\end{equation}
Using the above inequality together with inequality~(\ref{ineq_v}), we have:
\begin{equation}
  \langle y_2|u_2\rangle_T  \geq -\mathcal{V}(0) - \frac{\overline{\lambda}}{2\underline{\lambda}_B}\| y_2\|^2_{2,T}
    \label{P_PD_pass}
\end{equation}
From Lemma~\ref{theorem1} and in the view of closed-loop interconnection in Fig~\ref{feed_trans}, we know the lack of passivity in output of one system can be compensated by passivity in input of the other system. In other word, if~ $\delta_2+\epsilon_1>0$, then all results in Theorem~\ref{prop_main} are valid. Therefore, from equations~(\ref{virt-ENV}) and (\ref{P_PD_pass}), the passivity is guaranteed if:
\begin{equation}
   \underline{\lambda}_B-\frac{\overline{\lambda}}{2\underline{\lambda}_B}>0\Rightarrow \overline{\lambda} < 2\underline{\lambda}^2_B
    \label{P_PD_pass2}
\end{equation}
\noindent 3) In the case of ID, using the fact that an indefinite matrix has at least one positive eigenvalue and at least
one negative eigenvalue, we can assume that $\overline{\lambda}>0$ and all results for PD case are valid for ID case. 

In conclusion, the gain matrices can vary as long as inequality~(\ref{P_PD_pass2}) is satisfied. If the change, i.e. sign of $\dot{K}$ and $\dot{B}$, leads to a ND value for $\bar{P}$, inequality~(\ref{P_PD_pass2}) is automatically satisfied, and the passivity is preserved. If the change leads to PD value for $\bar{P}$ and inequality~(\ref{P_PD_pass2}) is effective, system stays passive. If none of above happens and change leads to an ID $\bar{P}$, having  inequality~(\ref{P_PD_pass2}) would be enough for the passivity. Note that the variation of damping matrix does not affect the passivity of $H_1$ in the interconnection of Fig.~\ref{feed_trans}(b). 

\end{proof}
\end{theorem}

\section{Concluding remarks and Future works}\label{dis_future_work}
In this research, a framework is established for controlling cooperative robots and executing the motion task of moving a load along a desired trajectory. A comprehensive model of an  augmented dynamics of the robot, JMs and the motors, provides the opportunity to introduce a new control scheme for semi-active joints called virtual control strategy (SVC). Based on SVC, any suitable control approach can be utilized to control the motion task in CRM. Here we used the concept of impedance control to devise the control scheme. The controller was studied and developed for two cases; constant damping and stiffness gains in the impedance relation; and variable gains. For both cases, the input/output passivity tool was used to analysis the stability.

Moreover, an optimization was introduced to obtain an energy-oriented impedance control. The optimization finds the best impedance gains such that the energy extraction from the power source is minimized. Using a simulation example, it was shown that energy regeneration can occur during the motion task, and it has huge effect in terms of energy saving. 

The lab experiment of this research is left as the future work. Also, another possible future work is to find an energy-optimal path for moving the load from point to point, which requires forming an optimization to find the optimal trajectory. 

\appendices
\section{ EXTERNAL ENERGY BALANCE IN CRM}\label{app_energy}
We start by obtaining Joule losses in terms of the desired control. The Joule losses due to resistance in motor of ij\textit{th} semi-active joint is:
\begin{equation}
L_{R_{ij}}=R_{ij}I_{{ij}}^2
\label{jl}
\end{equation}

Substituting current and $u_{ij}$ from equations~(\ref{Viq}) and (\ref{SVC}), respectively, and after some manipulation and simplification, we get:
\begin{equation}
L_{R_{ij}}=\frac{R_{ij}}{a_{ij}^2}(\tau^d_{ij})^2+\frac{a_{ij}^2\dot{q}_{ij}^2}{R_{ij}}-2\tau^d_{ij}\dot{q}_{ij}
\label{jl2}
\end{equation}

Equation~(\ref{JM_aug}) expresses the relation between $\tau_{ij}$ and $\tau^d_{ij}$. Multiplying both sides of this equation by $\dot{q}_{ij}$ yields:
\begin{equation}
\tau_{ij}\dot{q}_{ij}=-m_{ij}\bar{n}^2_{ij}\ddot{q}_{ij}\dot{q}_{ij}-(b_{ij}\bar{n}^2_{ij}+\frac{a_{ij}^2}{R_{ij}})\dot{q}^2_{ij}+\tau^d_{ij}\dot{q}_{ij}
\label{jl3}
\end{equation}

The kinetic energy of the actuator is expressed as $K_{ij}=\frac{1}{2}I_{ij}\bar{n}^2_{ij}\dot{q}^2_{ij}$. So equation~(\ref{jl3}) can be simplified as:
\begin{equation}
\tau_{ij}\dot{q}_{ij}=-\frac{dK_{ij}}{dt}-(b_{ij}\bar{n}^2_{ij}+\frac{a_{ij}^2}{R_{ij}})\dot{q}^2_{ij}+\tau^d_{ij}\dot{q}_{ij}
\label{jl4}
\end{equation}

Replacing $a_{ij}^2\dot{q}_{ij}/R_{ij}$ from equation~(\ref{jl3}) and rearranging the result and taking the integral from $t_1$ to $t_2$ of both sides of obtained equation gives:
\begin{equation}
\begin{split}
&\int_{t_1}^{t_2}(\tau^d_{ij}\dot{q}_{ij}-\frac{R_{ij}}{a_{ij}^2}(\tau^d_{ij})^2)dt= \\&\int_{t_1}^{t_2}(-\frac{dK_{ij}}{dt}-\tau_{ij}\dot{q}_{ij}-L_{R_{ij}}-b_{ij}\bar{n}^2_{ij}\dot{q}^2_{ij})dt
\end{split}
\label{jl5}
\end{equation}

On the other hand, using $\Delta E_c$ from equation~(\ref{internal_e_balance}), we can write:

\begin{equation}
\Delta E_c = \sum_{i=1}^{N}\sum_{j=1}^{n}\Delta E_{c_{ij}}
\label{jl6}
\end{equation}
where:
\begin{equation}
\Delta E_{c_{ij}}=\int_{t_1}^{t_2}(\tau^d_{ij}\dot{q}_{ij}-\frac{R_{ij}}{a_{ij}^2}(\tau^d_{ij})^2)dt
\label{jl7}
\end{equation}

Using equations~(\ref{jl5}) and (\ref{jl7}), it can be concluded:
\begin{equation}
\Delta E_{c_{ij}}=\int_{t_1}^{t_2}(-\frac{dK_{ij}}{dt}-\tau_{ij}\dot{q}_{ij}-L_{R_{ij}}-b_{ij}\bar{n}^2\dot{q}^2_{ij})dt
\label{jl8}
\end{equation}

To derive equation~(\ref{external_balance}), we start by writing the overall energy balance for CRM:
\begin{equation}
\int_{t_1}^{t_2} \dot{q}^T\mathcal{T}dt+\int_{t_1}^{t_2}\dot{q}^T\mathcal{T}_{ext}dt=\Delta E^\circ_m+\Sigma^\circ_m
\label{jl9}
\end{equation}
where the first term on the left-hand side is the work done by the semi-active joints, the second term is the work done by the external forces and moments ($W_{ext}$), $\Delta E^\circ_m$ is the total change in mechanical energy, and $\Sigma^\circ_m$ is the dissipated mechanical energy in the system. Also $\mathcal{T}$ is defined as the vector of all applied forces/moments and $\mathcal{T}_{ext}=J^T(q)f$ is the vector of external forces/moments, i.e. the forces/moments applied by the load.  

If we write equation~(\ref{jl5}) for all joints in CRM and substitute in equation~(\ref{jl9}), results:

\begin{equation}
W_{{ext}} = \Delta\,E_{c}\,+\,\Delta\,E^{Tot}_{m}\,+\,\Sigma^{Tot}_{m}\,+\,\Sigma_{e}
\label{jl10}
\end{equation}
where
\begin{equation}
\begin{split}
    \Delta\,E^{Tot}_{m} &= \Delta E^\circ_m+\sum_{i=1}^{N}\sum_{j=1}^{n}\Delta K_{ij} \\
    \Sigma^{Tot}_{m} &= \Sigma^\circ_m +\sum_{i=1}^{N}\sum_{j=1}^{n}\int_{t_1}^{t_2} b_{ij}\bar{n}^2_{ij}\dot{q}^2_{ij}dt\\ 
    \Sigma_e &=\sum_{i=1}^{N}\sum_{j=1}^{n} \int_{t_1}^{t_2}L_{R_{ij}}dt 
\end{split}
\label{jl11}
\end{equation}


\section*{Acknowledgment}
The authors would like to thank the National Science Foundation for funding this work (NSF grant $\#$1536035).

\ifCLASSOPTIONcaptionsoff
  \newpage
\fi

\end{document}